# A Note on Sampled-Data Observers

**Iasson Karafyllis[*], Tarek Ahmed-Ali[**] and Fouad Giri[**]**

[*]Dept. of Mathematics, National Technical University of Athens,
Zografou Campus, 15780, Athens, Greece,
email: iasonkar@central.ntua.gr; iasonkaraf@gmail.com

[**]Normandie UNIV, UNICAEN, ENSICAEN, LAC, 14000 Caen, France
email: tarek.ahmed-ali@ensicaen.fr; fouad.giri@unicaen.fr

**Abstract**

We present a new approach for deriving sampled-data observers from continuous-time observers that feature an Input-to-Output Stability property with respect to the output measurement noise and exponential convergence in the noiseless case. The design approach applies to a wide class of systems and yields sampled-data observers that inherit all performance characteristics of the underlying continuous-time observers. The main component of the proposed sampled-data observer is a novel output predictor that encompasses both inter-sample predictors and ZOH-predictors.

## 1. Introduction

Digital technology has nowadays gained almost all sectors of industry and society. This feature has affected control systems which have now become of hybrid (continuous-discrete) nature. Consequently, new research topics have become central in the theory of control systems. Such an example is the topic of compensation of time-delay and data-sampling effects in controller and observer design for continuous-time nonlinear systems. In this paper, we focus on the problem of designing sampled-data observers for nonlinear (continuous-time) systems. As a matter of fact, data-sampling entails data loss and introduces a varying time-delay. If not satisfactorily compensated for, these combined effects may result in the loss of system observability and divergence of the observer error. Several approaches have been developed to cope with this issue. One of them consists in letting the observer design be based on an Euler-like discrete-time approximation of the (continuous-time) system. The discrete-time observer obtained in this way only provides state estimates at the sampling instants. This approach has been investigated in [3,12] where it was shown that the corresponding observers ensure semi-global practical stability of the observation error. Observers that feature global exponential stability are generally designed directly from the (continuous-time) system model. The dominant design principle consists in starting from a (continuous-time) observer featuring global exponential stability, when continuous output measurements are available, and modifying the observer to account for data-sampling. The various modifications proposed in the literature amount to compensate for the missing information (between two successive sampling times) by using output and/or state inter-sample estimators. The simplest variant of this general principle consists in simply using a Zero-Order-Hold (ZOH) output predictor i.e. using the most recent output measurement until the acquisition of a new measurement. This approach has been investigated in [17] considering a class of systems with Lipschitz state nonlinearity. Sufficient conditions for a Luenberger-type observer to be globally exponentially convergent have been expressed in the form of Linear Matrix Inequalities

(LMIs) involving, particularly, the observer gain and the sampling period. To enlarge the set of admissible sampling periods, it was suggested to modify the standard ZOH-output-predictor based observer design by letting the observer gain be exponentially decaying within the sampling intervals and reset it at sampling times; see [1]. The idea of using an exponentially decaying factor in the observer gain has been first introduced in [4] where the problem of observer design has been investigated for output-delayed systems (not subject to output sampling).

A less simple approach is that commonly referred-to as continuous-discrete observer design, initiated in the early 90s in [5] and developed later in [15,16]. The proposed discrete-continuous-time observers have been derived either from the (continuous-time) high-gain observer or the extended Kalman filter. They involve open-loop continuous-time inter-sample state predictors and a discrete-time feedback correction of the state estimate operated at sampling times. The feedback correction of the state trajectory is performed by adding, in the observer state equation, an innovative term proportional to the output estimation error (between the system output and the observer output) that is amplified with the observer gain. In the case of Kalman filter like observer the gain is continuously updated. The exponential convergence of the observation output error is established under ad hoc assumptions depending on the underlying (continuous-time) observer.

A quite different approach, initiated in [8], provides sampled-data based Luenburger-type continuous-time observers, where the sampling effect is compensated for using inter-sample output predictors. In this approach, the output predictor is the only observer component that is reset at the sampling times and the observer state equation is continuously driven with an innovative term proportional to the (predicted) estimation error. Compared to continuous-discrete observers, output-predictor based observers feature global exponential convergence as well as implementation simplicity, since only one equation of the observer is reinitialized at the sampling times (namely, that of the output predictor) and since uncertain sampling schedules are allowed. This approach has been proved to be applicable to several classes of systems including linear detectable systems and triangular globally Lipschitz systems; see [8]. The inter-sample output-predictor based observer design has been the subject of several extensions; see [2,6,10,11,13,14]. It has been resorted to design sampled-data observers for systems with parameter uncertainty in [6] and for systems with time-delays in [2,7]. The inter-sample output-predictor principle has also been proved to be useful in designing observers for parabolic PDEs (see [11]) and it has also been used in output feedback design (see [10]). Its extension to the case of asynchronous measurements has been reported in [13,14].

In the light of the above review, it appears that the problem of designing sampled-data observers has mainly been addressed for two specific classes of (finite-dimensional) systems, specifically strict-feedback systems or state-affine systems. Furthermore, all proposed sampled-data observers are derived from previously existing exponentially-convergent continuous-time observers. The latter are made sampled-data using either the ZOH technique or the inter-sample output-predictor technique. In this paper, we revisit the sampled-data observer design problem with the aim of achieving a more general design approach that unifies most existing design methods and, in some aspects, goes much beyond the current results. The approach we present features the following facts:

(i) It applies to a large class of (locally) Lipschitz nonlinear systems with not-necessarily strict-feedback structure.
(ii) It involves a quite general observer structure that accounts for many known continuous-time observers. It is only required that the considered observer structure results in a (continuous-time) state error system that is Input-to-Output Stable (IOS) with respect to measurement noise and exponentially stable in the noiseless case.
(iii) It introduces a novel output predictor that encompasses both ZOH and the inter-sample predictors. The novel output predictor constitutes the main instrument in the sample-data observer and allows uncertain sampling schedules.
(iv) We provide an explicit condition on the maximum allowable sampling period so that the IOS property is preserved by the sampled-data observer.

The paper is organized as follows: in Section 2, we describe the main results of the present work as well as the design procedure of the proposed sampled-data observers. A worked example is considered in Section 3 that illustrates the supremacy of the inter-sample predictor design over other designs. All technical proofs are provided in Section 4. Some concluding remarks end the paper.

**Notation:** Throughout the paper, we adopt the following notation.

* $\Re_+ := [0,+\infty)$.

* Let $S \subseteq \Re^n$ be an open set and let $A \subseteq \Re^n$ be a set that satisfies $S \subseteq A \subseteq cl(S)$. By $C^0(A;\Omega)$, we denote the class of continuous functions on $A$, which take values in $\Omega \subseteq \Re^m$. By $C^k(A;\Omega)$, where $k \geq 1$ is an integer, we denote the class of functions on $A \subseteq \Re^n$, which takes values in $\Omega \subseteq \Re^m$ and has continuous derivatives of order $k$. In other words, the functions of class $C^k(A;\Omega)$ are the functions which have continuous derivatives of order $k$ in $S = \text{int}(A)$ that can be continued continuously to all points in $\partial S \cap A$. When $\Omega = \Re$ then we write $C^0(A)$ or $C^k(A)$. For $f \in C^0([0,1])$ the sup norm is defined by $\|f\|_\infty = \sup_{0 \leq z \leq 1} (|f(z)|) < +\infty$.

* For a vector $x \in \Re^n$ we denote by $|x|$ its usual Euclidean norm and by $x'$ its transpose. By $|A| := \sup\{|Ax|; x \in \Re^n, |x|=1\}$ we denote the induced norm of a matrix $A \in \Re^{m \times n}$ and $I$ denotes the identity matrix.

* We say that an increasing and continuous function $\rho : \Re_+ \to \Re_+$ is of class $K_\infty$ if $\rho(0) = 0$ and $\lim_{s \to +\infty} \rho(s) = +\infty$.

* Let $D \subseteq \Re^l$ be a non-empty set and $A \subseteq \Re_+$ an interval. By $L^\infty_{loc}(A;D)$ we denote the class of all Lebesgue measurable and locally bounded mappings $d : A \to D$. Notice that by $\sup_{\tau \in A}(|d(\tau)|)$ we do not mean the essential supremum of $d$ on $A$ but the actual supremum of $d$ on $A$.

## 2. Main Results

We consider forward complete systems of the form

$$\dot{x} = f(x,u) \qquad x \in \Re^n, u \in U \tag{2.1}$$

where $x \in \Re^n$ is the state of the system, $u \in U$ is a measured input, $U \subseteq \Re^m$ is a non-empty set and $f : \Re^n \times U \to \Re^n$ is a continuous mapping that satisfies a local Lipschitz condition, i.e., for every bounded set $S \subseteq \Re^n \times U$ there exists a constant $M(S) > 0$ such that $|f(x,u) - f(z,u)| \leq M(S)|x-z|$ for all $(x,u) \in S$, $(z,u) \in S$. By forward complete, we mean that for every $x_0 \in \Re^n$, $u \in L^\infty_{loc}(\Re_+;U)$ the unique solution $x(t) \in \Re^n$ of (2.1) with initial condition $x(0) = x_0$ corresponding to input $u \in L^\infty_{loc}(\Re_+;U)$ exists for all $t \geq 0$.

The measured output of system (2.1) is given by the equation

$$y = h(x) + \xi \qquad y \in \Re^p, \xi \in \Re^p \tag{2.2}$$

where $h: \Re^n \to \Re^p$ is a smooth mapping and $\xi$ is the measurement error (noise). We assume that when the output is measured continuously we can construct a continuous-time observer of the form

$$\dot{z} = f(z,u) + g(z,y,u)(y - h(z)) \qquad (2.3)$$
$$z \in \Re^n$$

where $z \in \Re^n$ is the state estimate and $g: \Re^n \times \Re^p \times U \to \Re^{n \times p}$ is a continuous mapping that satisfies a local Lipschitz condition, i.e., for every bounded set $S \subseteq \Re^n \times \Re^p \times U$ there exists a constant $M(S) > 0$ such that $|g(x,y,u) - g(z,w,u)| \leq M(S)(|x-z| + |y-w|)$ for all $(x,y,u) \in S$, $(z,w,u) \in S$.

When the output measurement is sampled, i.e., when the output values are available only at certain times which form an increasing sequence $\{t_k \geq 0: k = 0,1,2,...\}$ with $t_0 = 0$, $\lim_{k \to +\infty}(t_k) = +\infty$, with

$$y(t_k) = h(x(t_k)) + \xi(t_k), \; k = 0,1,2,... \qquad (2.4)$$

then the observer (2.3) has to be modified. We consider next three such modifications that have been proposed in the literature:

1) Observer with ZOH; see [17]. For this sampled-data observer, we replace the non-available signal $y(t) - h(z(t))$ by the most recently measured value of the signal, i.e., by $y(t_k) - h(z(t_k))$ for $t \in [t_k, t_{k+1})$. Equivalently, we replace the non-available signal $y(t)$ by $y(t_k) - h(z(t_k)) + h(z(t))$ for $t \in [t_k, t_{k+1})$ and the ZOH version of (2.3) is given by the equations:

$$\dot{z}(t) = f(z(t), u(t)) + g(z(t), w(t), u(t))(w(t) - h(z(t))) \qquad (2.5)$$

$$\dot{w}(t) = \nabla h(z(t))\dot{z}(t), t \in [t_k, t_{k+1}), \; k = 0,1,2,... \qquad (2.6)$$

$$w(t_k) = y(t_k), \; k = 0,1,2,... \qquad (2.7)$$

where $w(t) \in \Re^p$.

2) Observer with ZOH and exponentially time-varying gain; see [1]. For this sampled-data observer, we replace the non-available signal $y(t) - h(z(t))$ by an exponentially weighted signal which is based on the most recently measured value of the signal, i.e., by $\exp(-\eta(t - t_k))(y(t_k) - h(z(t_k)))$ for $t \in [t_k, t_{k+1})$, where $\eta \geq 0$ is a constant. Equivalently, we replace the non-available signal $y(t)$ by the signal $\exp(-\eta(t - t_k))(y(t_k) - h(z(t_k))) + h(z(t))$ for $t \in [t_k, t_{k+1})$ the ZOH version of (2.3) with exponentially time-varying gain will be given by the equations (2.5), (2.7) and

$$\dot{w}(t) = \nabla h(z(t))\dot{z}(t) - \eta(w(t) - h(z(t))), t \in [t_k, t_{k+1}), \; k = 0,1,2,... \qquad (2.8)$$

3) Observer with inter-sample predictor; see [8]. For this sampled-data observer, the non-available signal $y(t)$ is replaced by its approximation $w(t)$ given by (2.5), (2.7) and

$$\dot{w}(t) = \nabla h(z(t)) f(z(t), u(t)), t \in [t_k, t_{k+1}), \; k = 0,1,2,... \qquad (2.9)$$

It is clear that all sampled-data observers above are special cases of the sampled-data observer given by (2.5), (2.7) and

$$\dot{w}(t) = \nabla h(z(t))f(z(t),u(t)) - K(z(t),w(t),u(t))\big(w(t)-h(z(t))\big), t \in [t_k, t_{k+1}), \; k=0,1,2,\ldots \quad (2.10)$$

where $K(z,w,u) \in \Re^{p \times p}$ is a matrix with continuous entries that satisfies a local Lipschitz condition, i.e., for every bounded set $S \subseteq \Re^n \times \Re^p \times U$ there exists a constant $M(S) > 0$ such that $|K(x,y,u) - K(z,w,u)| \leq M(S)(|x-z| + |y-w|)$ for all $(x,y,u) \in S$, $(z,w,u) \in S$.

More specifically, the ZOH version of (2.3) is given by (2.5), (2.7), (2.10) with $K(z,w,u) = -\nabla h(z)g(z,w,u)$, the ZOH version of (2.3) with exponentially time-varying gain is given by (2.5), (2.7), (2.10) with $K(z,w,u) = -\nabla h(z)g(z,w,u) + \eta I$ and the sampled-data version of (2.3) with inter-sample predictor is given by (2.5), (2.7), (2.10) with $K(z,w,u) \equiv 0$.

Therefore, we need to study the sampled-data observer (2.5), (2.7), (2.10):
i) for theoretical reasons, since this observer unifies many proposed sampled-data observers, and
ii) for practical reasons, because an appropriate selection of the matrix $K(z,w,u) \in \Re^{p \times p}$ may allow larger sampling periods or reduced sensitivity to measurement noise.

To this purpose, we introduce the following technical assumption.

**(H)** *There exist continuous functions $V,W : \Re^n \times \Re^n \to \Re_+$, $a \in K_\infty$ and constants $\omega > 0$, $\gamma, L \geq 0$, $q \in \Re$ such that the following inequalities hold for all $z, x \in \Re^n$, $w \in \Re^p, u \in U$*

$$a(|z-x|) \leq V(z,x) \quad (2.11)$$

$$\begin{aligned}&(w-h(x))'\big(\nabla h(z)f(z,u) - \nabla h(x)f(x,u) - K(z,w,u)(w-h(z))\big) \\ &\leq LV(z,x) + q|w-h(x)|^2\end{aligned} \quad (2.12)$$

*and such that for every $x_0 \in \Re^n$, $z_0 \in \Re^n$, $u \in L^\infty_{loc}(\Re_+; U)$, $w \in L^\infty_{loc}(\Re_+; \Re^p)$ the unique solution of*

$$\begin{aligned}\dot{x} &= f(x,u) \\ \dot{z} &= f(z,u) + g(z,w,u)(w-h(z))\end{aligned} \quad (2.13)$$

*with initial condition $x(0) = x_0, z(0) = z_0$ corresponding to inputs $u \in L^\infty_{loc}(\Re_+; U)$, $w \in L^\infty_{loc}(\Re_+; \Re^p)$ exists for all $t \geq 0$ and satisfies the estimate:*

$$V(z(t),x(t)) \leq \exp(-\omega t)W(z(0),x(0)) + \gamma \sup_{0 \leq s \leq t}\big(|w(s)-h(x(s))|^2 \exp(-\omega(t-s))\big), \text{ for all } t \geq 0 \quad (2.14)$$

Assumption (H) is a generalized version of similar assumptions that have been used in [8]. Notice that estimate (2.14) essentially is an Input-to-Output Stability (IOS) property with respect to the measurement noise.

Our first main result is given next.

**Theorem 2.1:** *Consider system (2.1) and suppose that assumption (H) holds. Let $T > 0$ be a constant that satisfies*

$$2\gamma L \int_0^T \exp(2qs) ds < 1 \qquad (2.15)$$

*Then for every constant $\sigma \in (0, \omega]$ with $2\gamma L \int_0^T \exp((2q+\sigma)s) ds < 1$, for every $x_0 \in \Re^n$, $z_0 \in \Re^n$, $u \in L^\infty_{loc}(\Re_+; U)$, $\xi \in L^\infty_{loc}(\Re_+; \Re^p)$ and for every increasing sequence $\{t_k \geq 0 : k = 0, 1, 2, ...\}$ with $t_0 = 0$, $\lim_{k \to +\infty}(t_k) = +\infty$, $\sup\{t_{k+1} - t_k : k = 0, 1, 2, ...\} \leq T$, the unique solution $(x(t), z(t), w(t)) \in \Re^n \times \Re^n \times \Re^p$ of (2.1), (2.4), (2.5), (2.7), (2.10) corresponding to inputs $u \in L^\infty_{loc}(\Re_+; U)$, $\xi \in L^\infty_{loc}(\Re_+; \Re^p)$ satisfies the estimate:*

$$V(z(t), x(t)) \leq \Omega W(z(0), x(0)) \exp(-\sigma t) + \gamma \Omega \exp(\max(0, 2qT)) \sup_{0 \leq s \leq t} \left( |\xi(s)|^2 \right), \text{ for all } t \geq 0 \qquad (2.16)$$

*where $\Omega := \left( 1 - 2\gamma L \int_0^T \exp((2q+\sigma)s) ds \right)^{-1}$.*

It should be noticed that continuity of the function $j(\sigma) := 2\gamma L \int_0^T \exp((2q+\sigma)s) ds$ in conjunction with (2.15) guarantees the existence of sufficiently small $\sigma \in (0, \omega]$ such that $2\gamma L \int_0^T \exp((2q+\sigma)s) ds < 1$. The constant $\sigma \in (0, \omega]$ determines the convergence rate of the observer error to zero when measurement noise is absent. Finally, it should be noticed that the gain coefficient of the measurement noise in (2.16) is always higher than the gain coefficient $\gamma$ appearing in (2.14) for the continuous-time observer: the sampled-data observer is always more sensitive to noise than the continuous-time observer due to the limited flow of information from the system to the observer.

When the function $V$ involved in Assumption (H) is a quadratic function and the output map $h$ is linear (case of linear systems or globally Lipschitz systems) then we are in a position to provide more explicit and less conservative results.

**Theorem 2.2:** *Consider system (2.1) and suppose that there exist matrices $C \in \Re^{p \times n}$, $R \in \Re^{n \times p}$ with $R \neq 0$, $P \in \Re^{n \times n}$ being symmetric and positive definite and constants $\omega, L > 0$, $q \in \Re$ such that the following hold:*

$$h(x) = Cx, \text{ for all } x \in \Re^n \qquad (2.17)$$

$$e'P(f(x+e, u) - f(x, u) - RCe) \leq -\omega e'Pe, \text{ for all } x, e \in \Re^n \qquad (2.18)$$

$$|Cf(x+e, u) - Cf(x, u) - qCe|^2 \leq L^2 e'Pe, \text{ for all } x, e \in \Re^n \qquad (2.19)$$

*Let $T > 0$ be a constant that satisfies*

$$T < \frac{1}{q}\ln\left(1 + \frac{\omega q}{L\sqrt{|R'PR|}}\right), \text{ if } q > 0 \text{ or } -\omega^{-1}L\sqrt{|R'PR|} < q < 0 \quad (2.20)$$

$$T < \frac{\omega}{L\sqrt{|R'PR|}}, \text{ if } q = 0 \quad (2.21)$$

*Then there exist constants $\sigma, \Omega, \gamma > 0$ such that for every $x_0 \in \Re^n$, $z_0 \in \Re^n$, $u \in L_{loc}^\infty(\Re_+; U)$, $\xi \in L_{loc}^\infty(\Re_+; \Re^p)$ and for every increasing sequence $\{t_k \geq 0 : k = 0,1,2,...\}$ with $t_0 = 0$, $\lim_{k \to +\infty}(t_k) = +\infty$, $\sup\{t_{k+1} - t_k : k = 0,1,2,...\} \leq T$, the unique solution of (2.1), (2.4), (2.5), (2.7), (2.10) with $K(z,w,u) \equiv -qI$, $g(z,w,u) \equiv R$, corresponding to inputs $u \in L_{loc}^\infty(\Re_+; U)$, $\xi \in L_{loc}^\infty(\Re_+; \Re^p)$ satisfies the estimate:*

$$|z(t) - x(t)| \leq \Omega \exp(-\sigma t)|z_0 - x_0| + \gamma \sup_{0 \leq s \leq t}(|\xi(s)|), \text{ for all } t \geq 0 \quad (2.22)$$

It should be noticed that there is no restriction on the values of the diameter of the sampling schedule $T > 0$ when $q \leq -\omega^{-1}L\sqrt{|R'PR|}$.

## 3. Illustrative Examples

As mentioned in the Introduction, the proposed sampled-data observer design can allow larger values for the sampling period when the matrix $K$ is selected in an appropriate way. To see this, consider the linear system

$$\begin{aligned} \dot{x}_1 &= x_2 \\ \dot{x}_2 &= -x_1 + u \\ y &= x_1 \end{aligned} \quad (3.1)$$

This linear system was selected because there is a natural upper bound of the upper diameter of the sampling schedule: the system is not observable for uniform sampling schedules $t_k = k\pi$, $k = 0,1,2,...$ and this implies that for every possible sampled-data observer design we must necessarily have $\sup\{t_{k+1} - t_k : k = 0,1,2,...\} < \pi$. For this system all assumptions of Theorem 2.2 hold with

$$C = \begin{bmatrix} 1 & 0 \end{bmatrix}, \omega = 1$$

$$P = \frac{1}{2}\begin{bmatrix} 2 & -1 \\ -1 & 1 \end{bmatrix}, R = \begin{bmatrix} 2 \\ 1 \end{bmatrix}$$

$$L = \sqrt{2(1 + (1-q)^2)}$$

Therefore, we are in a position to select $q \in \Re$ in an optimal way so that the upper diameter $T > 0$ of the sampling schedule becomes as large as possible. To this purpose, we define the function:

$$T_{\max}(q) = \frac{1}{q}\ln\left(1 + \frac{q}{\sqrt{1+(1-q)^2}\sqrt{5}}\right), \text{ if } q \neq 0 \tag{3.2}$$

$$T_{\max}(0) = \frac{1}{\sqrt{10}}, \tag{3.3}$$

Notice that since $|R'PR| = \frac{5}{2}$, it follows from Theorem 2.2 that the upper diameter $T > 0$ of the sampling schedule must satisfy the inequality $T < T_{\max}(q)$. The behavior of the function $T_{\max}(q)$ is shown in Fig. 1, where it is shown that the function presents a global maximum at $q = 0.8$ with $T_{\max}(q) = 0.37589$. Therefore, when $q = 0.8$, the sampled-data observer

$$\begin{aligned}\dot{z}_1(t) &= -2z_1(t) + z_2(t) + 2w(t) \\ \dot{z}_2(t) &= -2z_1(t) + u(t) + w(t)\end{aligned} \tag{3.4}$$

$$\dot{w}(t) = z_2(t) + q(w(t) - z_1(t)), t \in [t_k, t_{k+1}), \ k = 0, 1, 2, \ldots \tag{3.5}$$

$$w(t_k) = x_1(t_k), \ k = 0, 1, 2, \ldots \tag{3.6}$$

will guarantee exponential convergence of the observer error $|z(t) - x(t)|$ provided that the sampling times form an increasing sequence $\{t_k \geq 0 : k = 0, 1, 2, \ldots\}$ with $t_0 = 0$, $\lim_{k \to +\infty}(t_k) = +\infty$, $\sup\{t_{k+1} - t_k : k = 0, 1, 2, \ldots\} < 0.37589$.

To understand the importance of the newly proposed sampled-data observer (3.4), (3.5), (3.6), it is useful to compare with two existing cases in the literature:

a) the ZOH version of the continuous-time observer

$$\begin{aligned}\dot{z}_1(t) &= -2z_1(t) + z_2(t) + 2x_1(t) \\ \dot{z}_2(t) &= -2z_1(t) + u(t) + x_1(t)\end{aligned} \tag{3.7}$$

which is the continuous-time observer with gains $R = \begin{bmatrix}2 & 1\end{bmatrix}'$ (used for the design of (3.4), (3.5), (3.6)), i.e., the sampled-data observer

$$\begin{aligned}\dot{z}_1(t) &= z_2(t) - 2(z_1(t_k) - x_1(t_k)) \\ \dot{z}_2(t) &= -z_1(t) + u(t) - (z_1(t_k) - x_1(t_k))\end{aligned} \tag{3.8}$$

which corresponds to the sampled-data observer (3.4), (3.5), (3.6) with $q = 2$. For $q = 2$ we have $T_{\max}(q) = 0.24504$ and it is clear that the sampled-data observer (3.4), (3.5), (3.6) with $q = 0.8$ allows a 53.4% increase of the diameter of the sampling schedule.

b) the continuous-time observer (3.7) with an inter-sample predictor, i.e., the sampled-data observer (3.4), (3.5), (3.6) with $q = 0$. For $q = 0$ we have $T_{\max}(q) = \frac{1}{\sqrt{10}} \approx 0.31623$ and it is clear that the sampled-data observer (3.4), (3.5), (3.6) with $q = 0.8$ allows a 18.8% increase of the diameter of the sampling schedule.

It should be noticed that for this simple example the sampled-data observer (3.4), (3.5), (3.6) is not different from the observer with ZOH and exponentially time-varying gain proposed in [1]. However, here we know that an exponential stability estimate for the observer error holds without having to solve LMIs. Moreover, using Theorem 2.2 we have selected optimally the exponential rate of the gain (i.e., $q$) in order to maximize the sampling period.

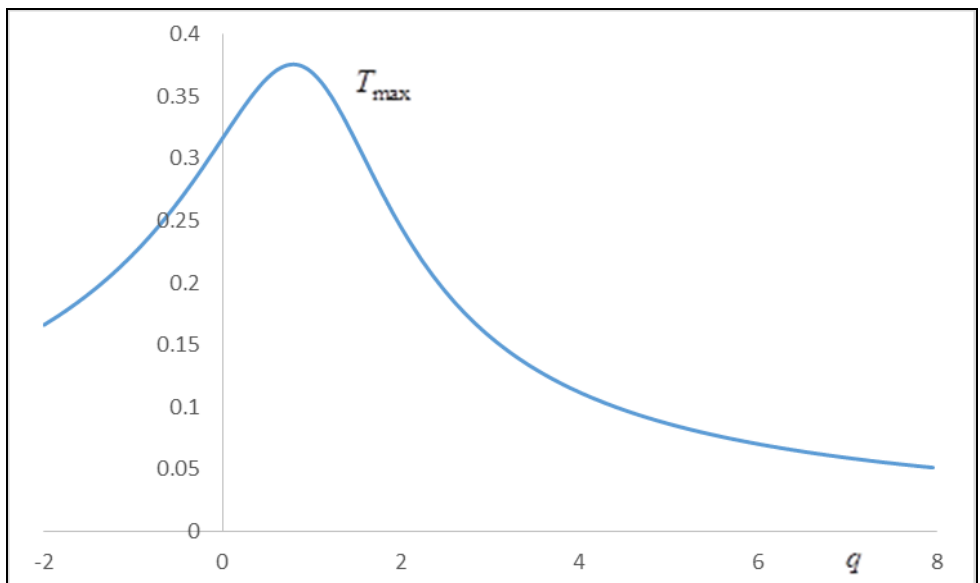

**Figure 1:** The graph of the function $T_{\max}(q)$ defined by (3.2), (3.3).

## 4. Proofs of Main Results

We next provide the proofs of Theorem 2.1 and Theorem 2.2.

**Proof of Theorem 2.1:** Let arbitrary $\sigma \in (0, \omega]$ with $2\gamma L \int_0^T \exp((2q+\sigma)s) ds < 1$, $x_0 \in \Re^n$, $z_0 \in \Re^n$, $u \in L^{\infty}_{loc}(\Re_+; U)$, $\xi \in L^{\infty}_{loc}(\Re_+; \Re^p)$ be given and let an arbitrary increasing sequence $\{t_k \geq 0 : k = 0,1,2,...\}$ with $t_0 = 0$, $\lim_{k \to +\infty}(t_k) = +\infty$, $\sup\{t_{k+1} - t_k : k = 0,1,2,...\} \leq T$ be also given. Consider an arbitrary sampling time $t_k$ and suppose that the solution $(x(t), z(t), w(t)) \in \Re^n \times \Re^n \times \Re^p$ of (2.1), (2.4), (2.5), (2.7), (2.10) corresponding to inputs $u \in L^{\infty}_{loc}(\Re_+; U)$, $\xi \in L^{\infty}_{loc}(\Re_+; \Re^p)$ exists on $[0, t_k]$. We will show next that the solution $(x(t), z(t), w(t)) \in \Re^n \times \Re^n \times \Re^p$ exists on $[0, t_{k+1}]$, which will imply that the solution exists for all $t \geq 0$.

Notice that the component $x(t)$ of the solution exists on $[0, t_{k+1}]$ (by virtue of forward completeness of (2.1)). Moreover, standard theory of ordinary equations guarantees that there exists $\delta > 0$ such that the solution of (2.5), (2.10) exists on $[t_k, t_k + \delta)$. Notice that inequality (2.12) implies the following differential inequality for $t \in [t_k, t_k + \delta)$ a.e.:

$$(w(t) - h(x(t)))' \frac{d}{dt}(w(t) - h(x(t))) \leq LV(z(t), x(t)) + q|w(t) - h(x(t))|^2 \quad (4.2)$$

Using Lemma 2.12 in [9] in conjunction with (4.2) and the fact that $w(t_k) - h(x(t_k)) = \xi(t_k)$ (a consequence of (2.4), (2.7)), we get for all $t \in [t_k, t_k + \delta)$:

$$|w(t) - h(x(t))|^2 \leq \exp(2q(t - t_k))|\xi(t_k)|^2 + 2L\int_{t_k}^{t} \exp(2q(t - s))V(z(s), x(s))ds \quad (4.3)$$

Estimate (4.3) implies the following estimate for all $t \in [t_k, t_k + \delta)$:

$$|w(t) - h(x(t))|^2 \exp(\sigma t) \leq \exp(2q(t - t_k))|\xi(t_k)|^2 \exp(\sigma t)$$
$$+ 2L\int_{t_k}^{t} \exp((2q + \sigma)(t - s))ds \sup_{t_k \leq s \leq t}(V(z(s), x(s))\exp(\sigma s)) \quad (4.4)$$

Using the fact that $\sup\{t_{k+1} - t_k : k = 0, 1, 2, ...\} \leq T$ in conjunction with estimate (4.4), we get for all $t \in [t_k, t_k + \delta)$:

$$|w(t) - h(x(t))|^2 \exp(\sigma t) \leq \exp(\max(0, 2qT))\exp(\sigma t) \sup_{0 \leq s \leq t}(|\xi(s)|^2)$$
$$+ 2L\int_{0}^{T} \exp((2q + \sigma)s)ds \sup_{0 \leq s \leq t}(V(z(s), x(s))\exp(\sigma s)) \quad (4.5)$$

Since estimate (4.5) is independent of the integer $k$, it follows that (4.5) holds for all $t \in [0, t_k + \delta)$.

Notice that the fact that $\sigma \in (0, \omega]$ and estimate (2.14) imply the estimate:

$$V(z(t), x(t))\exp(\sigma t) \leq W(z(0), x(0)) + \gamma \sup_{0 \leq s \leq t}(|w(s) - h(x(s))|^2 \exp(\sigma s)), \text{ for all } t \in [0, t_k + \delta) \quad (4.6)$$

Combining (4.5) and (4.6) we obtain for all $t \in [0, t_k + \delta)$:

$$V(z(t), x(t))\exp(\sigma t) \leq W(z(0), x(0)) + \gamma \exp(\max(0, 2qT))\exp(\sigma t) \sup_{0 \leq s \leq t}(|\xi(s)|^2)$$
$$+ 2\gamma L\int_{0}^{T} \exp((2q + \sigma)s)ds \sup_{0 \leq s \leq t}(V(z(s), x(s))\exp(\sigma s)) \quad (4.7)$$

Since $2\gamma L\int_{0}^{T} \exp((2q + \sigma)s)ds < 1$, it follows from (4.7) that the following estimate holds for all $t \in [0, t_k + \delta)$:

$$\sup_{0 \leq s \leq t}(V(z(s), x(s))\exp(\sigma s)) \leq \left(1 - 2\gamma L\int_{0}^{T} \exp((2q + \sigma)s)ds\right)^{-1} W(z(0), x(0))$$
$$+ \gamma\left(1 - 2\gamma L\int_{0}^{T} \exp((2q + \sigma)s)ds\right)^{-1} \exp(\max(0, 2qT))\exp(\sigma t) \sup_{0 \leq s \leq t}(|\xi(s)|^2) \quad (4.8)$$

Estimate (4.8) and the fact that $x(t)$ is bounded on $[0, t_{k+1}]$, in conjunction with (2.11) implies that $z(t)$ is bounded on $[0, t_k + \delta)$. Similarly, combining (4.8) and (4.5), we conclude that $w(t)$ is bounded on $[0, t_k + \delta)$. Therefore, the theory of ordinary differential equations allows us to conclude that the solution of (2.5), (2.10) exists on $[t_k, t_{k+1}]$. Replacing the value of $w(t_{k+1})$ by using (2.4), (2.7), we conclude that the solution $(x(t), z(t), w(t)) \in \Re^n \times \Re^n \times \Re^p$ of (2.1), (2.4), (2.5), (2.7), (2.10) corresponding to inputs $u \in L^\infty_{loc}(\Re_+; U)$, $\xi \in L^\infty_{loc}(\Re_+; \Re^p)$ exists on $[0, t_{k+1}]$.

We conclude that (4.8) holds for all $t \geq 0$. Therefore, estimate (2.16) holds for all $t \geq 0$. The proof is complete. ◁

**Proof of Theorem 2.2:** Notice that inequalities (2.20), (2.21) guarantee that

$$\frac{L^2 |R'PR|}{\omega^2} \left( \int_0^T \exp(qs) ds \right)^2 < 1 \tag{4.9}$$

By virtue of (4.9) and continuity of the function $j(\sigma) := \frac{L^2 |R'PR|}{\omega(\omega - 2\sigma)} \left( \int_0^T \exp((q+\sigma)s) ds \right)^2$ at $\sigma = 0$, there exists $\sigma \in (0, \omega/2)$ such that

$$\frac{L^2 |R'PR|}{\omega(\omega - 2\sigma)} \left( \int_0^T \exp((q+\sigma)s) ds \right)^2 < 1 \tag{4.10}$$

Let arbitrary $x_0 \in \Re^n$, $z_0 \in \Re^n$, $u \in L^\infty_{loc}(\Re_+; U)$, $\xi \in L^\infty_{loc}(\Re_+; \Re^p)$ be given and let an arbitrary increasing sequence $\{t_k \geq 0 : k = 0, 1, 2, ...\}$ with $t_0 = 0$, $\lim_{k \to +\infty}(t_k) = +\infty$, $\sup\{t_{k+1} - t_k : k = 0, 1, 2, ...\} \leq T$ be also given. Consider an arbitrary sampling time $t_k$ and suppose that the solution $(x(t), z(t), w(t)) \in \Re^n \times \Re^n \times \Re^p$ of (2.1), (2.4), (2.5), (2.7), (2.10) with $K(z, w, u) \equiv -qI$, $g(z, w, u) \equiv R$, corresponding to inputs $u \in L^\infty_{loc}(\Re_+; U)$, $\xi \in L^\infty_{loc}(\Re_+; \Re^p)$ exists on $[0, t_k]$. We will show next that the solution $(x(t), z(t), w(t)) \in \Re^n \times \Re^n \times \Re^p$ exists on $[0, t_{k+1}]$, which will imply that the solution exists for all $t \geq 0$.

Notice that the component $x(t)$ of the solution exists on $[0, t_{k+1}]$ (by virtue of forward completeness of (2.1)). Moreover, standard theory of ordinary equations guarantees that there exists $\delta > 0$ such that the solution of (2.5), (2.10) exists on $[t_k, t_k + \delta)$. Notice that equations (2.4), (2.7), (2.10), (2.17) and the fact that $K(z, w, u) \equiv -qI$, imply the following equation for all $t \in [t_k, t_k + \delta)$:

$$w(t) - Cx(t) = \exp(q(t - t_k)) \xi(t_k) \\ + \int_{t_k}^t \exp(q(t-s)) C \big( f(z(s), u(s)) - f(x(s), u(s)) + q(x(s) - z(s)) \big) ds \tag{4.11}$$

Using (2.19) and (4.11), we obtain the following estimate for all $t \in [t_k, t_k + \delta)$:

$$|w(t) - Cx(t)| \leq \exp(q(t - t_k)) |\xi(t_k)| + L \int_{t_k}^t \exp(q(t-s)) \sqrt{e'(s) P e(s)} ds \tag{4.12}$$

where $e(t) := z(t) - x(t)$. Estimate (4.12) implies the following estimate for all $t \in [t_k, t_k + \delta)$:

$$|w(t)-Cx(t)|\exp(\sigma t) \leq \exp(q(t-t_k))|\xi(t_k)|\exp(\sigma t)$$
$$+L\int_{t_k}^{t}\exp((q+\sigma)(t-s))ds \sup_{t_k\leq s\leq t}\left(\sqrt{e'(s)Pe(s)}\exp(\sigma s)\right) \quad (4.13)$$

Using the fact that $\sup\{t_{k+1}-t_k : k=0,1,2,...\} \leq T$ in conjunction with estimate (4.4), we get for all $t \in [t_k, t_k+\delta)$:

$$|w(t)-Cx(t)|\exp(\sigma t) \leq \exp(\max(0,qT)) \sup_{0\leq s\leq t}(|\xi(s)|)\exp(\sigma t)$$
$$+L\int_{0}^{T}\exp((q+\sigma)s)ds \sup_{0\leq s\leq t}\left(\sqrt{e'(s)Pe(s)}\exp(\sigma s)\right) \quad (4.14)$$

Since estimate (4.14) is independent of the integer $k$, it follows that (4.14) holds for all $t \in [0, t_k+\delta)$.

It follows from (2.1), (2.5), the facts that $e(t):=z(t)-x(t)$, $g(z,w,u) \equiv R$ and (2.17), (2.18) that the following differential inequality holds for $t \in [0, t_k+\delta)$ a.e.:

$$\frac{d}{dt}e'(t)Pe(t) \leq -2\omega e'(t)Pe(t) + 2e'(t)PR(w(t)-Cx(t)) \quad (4.15)$$

Using the fact that $\varepsilon a'Pa + \varepsilon^{-1}b'Pb \geq 2b'Pa$ for all $\varepsilon > 0$, $a,b \in \mathfrak{R}^n$, we obtain from (4.15) for $t \in [0, t_k+\delta)$ a.e. (with $\varepsilon = \omega > 0$, $a = e(t)$, $b = R(w(t)-Cx(t))$):

$$\frac{d}{dt}e'(t)Pe(t) \leq -\omega e'(t)Pe(t) + \omega^{-1}|R'PR||w(t)-Cx(t)|^2 \quad (4.16)$$

Using Lemma 2.12 in [9] in conjunction with (4.16), we get for all $t \in [0, t_k+\delta)$:

$$e'(t)Pe(t) \leq \exp(-\omega t)e'(0)Pe(0)$$
$$+\omega^{-1}|R'PR|\int_{0}^{t}\exp(-\omega(t-s))|w(s)-h(x(s))|^2 ds \quad (4.17)$$

Estimate (4.17) in conjunction with the fact that $\sigma \in (0, \omega/2)$ and the fact that $\sqrt{a+b} \leq \sqrt{a}+\sqrt{b}$ for all $a,b \geq 0$, give the following estimate for all $t \in [0, t_k+\delta)$:

$$\sqrt{e'(t)Pe(t)}\exp(\sigma t) \leq \sqrt{e'(0)Pe(0)} + \sqrt{\frac{|R'PR|}{\omega(\omega-2\sigma)}} \sup_{0\leq s\leq t}(|w(s)-h(x(s))|\exp(\sigma s)) \quad (4.18)$$

Combining estimates (4.14), (4.18), we obtain for all $t \in [0, t_k+\delta)$:

$$\sqrt{e'(t)Pe(t)}\exp(\sigma t) \leq \sqrt{e'(0)Pe(0)} + \exp(\max(0,qT))\sqrt{\frac{|R'PR|}{\omega(\omega-2\sigma)}} \sup_{0\leq s\leq t}(|\xi(s)|)\exp(\sigma t)$$
$$+L\sqrt{\frac{|R'PR|}{\omega(\omega-2\sigma)}}\int_{0}^{T}\exp((q+\sigma)s)ds \sup_{0\leq s\leq t}\left(\sqrt{e'(s)Pe(s)}\exp(\sigma s)\right) \quad (4.19)$$

Using (4.10) in conjunction with (4.19), we obtain for all $t \in [0, t_k+\delta)$:

$$\sup_{0\leq s\leq t}\left(\sqrt{e'(s)Pe(s)}\exp(\sigma s)\right)\leq \left(1-L\sqrt{\frac{|R'PR|}{\omega(\omega-2\sigma)}}\int_0^T \exp((q+\sigma)s)ds\right)^{-1}\sqrt{e'(0)Pe(0)}$$

$$+\exp(\max(0,qT)+\sigma t)\left(1-L\sqrt{\frac{|R'PR|}{\omega(\omega-2\sigma)}}\int_0^T \exp((q+\sigma)s)ds\right)^{-1}\sqrt{\frac{|R'PR|}{\omega(\omega-2\sigma)}}\sup_{0\leq s\leq t}\left(|\xi(s)|\right)$$

(4.20)

Estimate (4.20) and the fact that $x(t)$ is bounded on $[0,t_{k+1}]$, in conjunction with the fact that $P\in \Re^{n\times n}$ is positive definite, implies that $z(t)$ is bounded on $[0,t_k+\delta)$. Similarly, combining (4.14) and (4.20), we conclude that $w(t)$ is bounded on $[0,t_k+\delta)$. Therefore, the theory of ordinary differential equations allows us to conclude that the solution of (2.5), (2.10) with $K(z,w,u)\equiv -qI$, $g(z,w,u)\equiv R$, exists on $[t_k,t_{k+1}]$. Replacing the value of $w(t_{k+1})$ by using (2.4), (2.7), we conclude that the solution $(x(t),z(t),w(t))\in \Re^n \times \Re^n \times \Re^p$ of (2.1), (2.4), (2.5), (2.7), (2.10) with $K(z,w,u)\equiv -qI$, $g(z,w,u)\equiv R$, corresponding to inputs $u\in L_{loc}^\infty(\Re_+;U)$, $\xi\in L_{loc}^\infty(\Re_+;\Re^p)$ exists on $[0,t_{k+1}]$.

We conclude that (4.20) holds for all $t\geq 0$. Since $P\in \Re^{n\times n}$ is positive definite, there exist constants $0<c_1\leq c_2$ such that $c_1|e|^2 \leq e'Pe \leq c_2|e|^2$ for all $e\in \Re^n$. Combining the previous inequality with (4.20), we obtain estimate (2.22) for appropriate constants $\Omega,\gamma >0$. The proof is complete. ◁

## 5. Concluding Remarks

The present work provided a new approach for deriving sampled-data observers from continuous-time observers that feature an IOS property with respect to the output measurement noise and exponential convergence in the noiseless case. The proposed approach unifies many existing approaches for the design of sampled-data observers and applies to a wide class of systems. The corresponding sampled-data observers inherit all performance characteristics of the underlying continuous-time observers. Future work may address the problem of relaxing Assumption (H) and consequently allow the construction of sampled-data observers for even wider classes of nonlinear systems.